\documentclass[a4,10pt,epsfig]{article}
\usepackage{graphicx,epstopdf}

\usepackage[a4paper, total={5.25in, 8in}]{geometry}

\usepackage{epsfig,epsf,epstopdf}
\usepackage{relsize}
\usepackage[pdf]{pstricks}
\usepackage{float}
\usepackage{mathptmx}
\usepackage{amsmath}
\usepackage{amssymb}
\usepackage{amsfonts}
\usepackage{mathtools}
\graphicspath{ {plots/} }
\usepackage{euscript}
\usepackage{textcomp}
\usepackage[subtle]{savetrees}
\usepackage{tikz, pgfplots}
\usetikzlibrary{shapes,arrows}
\usepackage{fancyhdr}
\usetikzlibrary{shadows,calc}
\pgfmathdeclarefunction{gauss}{2}{%
  \pgfmathparse{1/(#2*sqrt(2*pi))*exp(-((x-#1)^2)/(2*#2^2))}%
}

\newcommand\qedsymbol{\hbox{\rlap{$\sqcap$}$\sqcup$}}
\newcommand\qed{\relax\ifmmode\else\unskip\hfill\fi{\Large $\sqbullet$}}
\newcommand\smartqed{\renewcommand\qed{\relax\ifmmode\qedsymbol\else
  {\unskip\nobreak\hfil\penalty50\hskip1em\null\nobreak\hfil\qedsymbol
  \parfillskip=\z@\finalhyphendemerits=0\endgraf}\fi}}

\abovedisplayshortskip
\belowdisplayshortskip
\abovedisplayskip
\belowdisplayskip

\newcommand{\bea}{\begin{eqnarray}}
\newcommand{\eea}{\end{eqnarray}}
\newcommand{\nbea}{\begin{eqnarray*}}
\newcommand{\neea}{\end{eqnarray*}}

\newtheorem{theorem}{Theorem}[section]
\newtheorem{corollary}{Corollary}[section]

\newtheorem{remark}{Remark}[section]

\begin{document}

\title{Decay estimate of bivariate Chebyshev coefficients for functions with limited smoothness}
\author{ Akansha\footnote{akansha@math.iitb.ac.in}\\
	Department of Mathematics, \\Indian Institute of Technology Bombay,\\ Powai, Mumbai - 400076. India.}
\date{}
\maketitle{}


\begin{abstract}
We obtain the decay bounds for Chebyshev series coefficients of functions with finite Vitali variation on the unit square. A generalization of the well known identity, which relates exact and approximated coefficients, obtained using the quadrature formula, is derived. Finally, an asymptotic $L^1$-approximation error of finite partial sum for functions of bounded variation in sense of Vitali as well as Hardy-Krause, on the unit square is deduced. 
\end{abstract}

\section{Introduction} 
The Chebyshev approximation is an optimal method to approximate smooth and non-smooth (Chebfun\footnote{Chebfun is an object oriented system written in MATLAB which extend MATLAB's basic commands, using Chebyshev approximation, to continuous context.} \cite{dri_hal_tre-14a,pac_pla_tre-10a}) functions. These trigonometric polynomial approximants are widely used in numerical schemes for solving partial differential equations. 

There are several results available for univariate functions which show that the decay rate of Chebyshev series coefficients depend on the smoothness of the function. For instance, we have the following decay bounds for 
functions in $W^{1,1}([-1,1])$ (\cite{tre-08a,xia_che_wan-10a,has-17a}) and $BV([-1,1])$ (\cite{xia-18a}).
\begin{enumerate}
\item \cite{tre-08a} If $f,f',\ldots,f^{(k-1)}$ is absolutely continuous on $[-1,1]$ and $V_k =\|f^{(k)}\|_T<\infty$ then for $j> k$
$$|c_j|\le \dfrac{2V_k}{\pi j(j-1)\cdots(j-k)},$$
where $k$ is a nonnegative integer, and $\|f\|_T$ is the weighted Chebyshev norm of $f$. 
The bounds obtained in \cite{tre-08a} are improved by Xiang \cite{xia_che_wan-10a}, and further, following sharper bounds are derived by Majidian \cite{has-17a}.
\item For a nonnegative integer $k$, if $f,f',\ldots,f^{(k-1)}$ is absolutely continuous on $[-1,1]$ and $V_k =\|f^{(k)}\|_T<\infty$ then for $j> k$
\begin{equation*}
|c_j| \le \frac{2V_k}{\pi}\begin{cases}
\dfrac{1}{j(j+2)(j-2)\cdots(j+2s)(j-2s)}, ~~~~~~~~~\mbox{if }k=2s,\\
\dfrac{1}{(j-1)(j+1)(j-3)\cdots(j+2s-1)(j-2s-1)}, ~~~~~\mbox{if }k=2s+1.
\end{cases}
\end{equation*}
\end{enumerate}
Xiang \cite{xia-18a} obtained the sharpest decay bounds for functions of limited regularities. So far, similar decay estimates for bivariate Chebyshev series coefficients of functions of bounded variation are not available in the literature. This article aims to extend the decay results obtained in \cite{has-17a} for univariate functions to two dimensions. The decay estimates for bivariate Fourier coefficients are obtained in \cite{gol_kal_lut-13a} for a class of smooth functions.  We derive decay estimates for bivariate Chebyshev coefficients of a large class of functions, including piecewise-smooth function with finite Vitali Variation. There are several generalizations of the definition of univariate bounded variation (for details, see \cite{cla_ada-33a,ada_cla-34a}); however, Vitali (later with additional conditions\footnote{The function $f(x,y_0)$ is of bounded variation for each $y_0\in[-1,1]$ and $f(x_0,y)$ is of bounded variation for each $x_0\in[-1,1]$} Hardy-Krause) proposed the most natural generalization, which 
preserves most of the properties of (our interest) univariate functions of bounded variation. An $L^1$ convergence result for functions of finite Vitali variation on the unit square is obtained. We also derive a relation between the exact and approximated bivariate Chebyshev coefficients to obtain the $L^1$-error estimate for Chebyshev approximants (with approximated coefficients) of functions of bounded variation in the sense of Hardy-Krause, on the unit square.

%

 
\section{Bivariate Chebyshev Series Expansion}\label{sec:bivariatechebyshev}
A continuous function $f$ of bounded variation with one of its partial derivative bounded on the unit square $D:=[-1,1]^2$ can be represented by the Chebyshev series as \cite{mas-han-03a} 

\begin{equation}\label{eq:bivariChebexpansion}
f(x,y)=\sum_{i = 0}^{\infty}{\vphantom{\sum}}'\sum_{j = 0}^{\infty}{\vphantom{\sum}}'c_{i,j}T_i(x)T_j(y),
\end{equation}
where $T_i(x)=\cos(i\cos^{-1}x)$ and $T_j(y)=\cos(j\cos^{-1}y)$ are the Chebyshev polynomial of degrees $i$ and $j$, respectively,
\begin{equation}\label{eq:bivarChebcoeffexact}
	c_{i,j} =\dfrac{4}{\pi^2} \int_Df(x,y)T_i(x)T_j(y)\omega(x,y)dxdy= \dfrac{4}{\pi^2}\int_{0}^{\pi}\int_{0}^{\pi}f(\cos\theta_x,\cos\theta_y)\cos i\theta_x \cos j\theta_y d\theta_x d\theta_y,
	\end{equation}
and	$\omega(x,y)=(1-x^2)^{-\frac{1}{2}}(1-y^2)^{-\frac{1}{2}}$ is the Chebyshev weight function. The prime in the summation indicates that coefficients of the terms $T_0(x)T_0(y)$, $T_0(x)T_j(y)$, for $j>0$, and $T_i(x)T_0(y)$, for all $i>0$, are $c_{0,0}/4$, $c_{0,j/2},$ and $c_{i,0}/2$, respectively.

Since the exact value of the integral \eqref{eq:bivarChebcoeffexact} for an arbitrary function cannot be obtained always, we use the following Gauss-Chebyshev quadrature rule 
\begin{equation}\label{eq:bivarChebcoeffapprx}
c_{i,j} \approx \dfrac{4}{n_xn_y}\sum_{l_x=1}^{n_x}\sum_{l_y=1}^{n_y}f(x_{l_x},y_{l_y})T_i(x_{l_x})T_j(y_{l_y}) =: \tilde{c}_{i,j},
\end{equation}
to approximate the bivariate Chebyshev series coefficients, where $x_{l_x}$ and $y_{l_y}$ are the roots of the Chebyshev polynomials $T_{n_x}(x)$ and $T_{n_y}(y)$, respectively. We denote a polynomial approximant of $f$ on $D$ with approximated and exact Chebyshev coefficients by 
\begin{equation*}\label{eq:bivarChebapprxexpansion}
\mathsf{\tilde{C}}_{d_x,d_y}[f](x,y) := \sum_{i = 0}^{d_x}{\vphantom{\sum}}'\sum_{j = 0}^{d_y}{\vphantom{\sum}}'\tilde{c}_{i,j}T_i(x)T_j(y),~
\mbox{ and }
~\mathsf{{C}}_{d_x,d_y}[f](x,y) := \sum_{i = 0}^{d_x}{\vphantom{\sum}}'\sum_{j = 0}^{d_y}{\vphantom{\sum}}'{c}_{i,j}T_i(x)T_j(y).
\end{equation*}
respectively. Note that, a bivariate Chebyshev approximant for a function defined on an arbitrary rectangular domain can be obtained using the standard change of variable process.
\section{Decay bounds and Convergence results}\label{sec:decay.bounds2Dche}
This section first obtains the decay bounds for Chebyshev series coefficients of a function $f:\Omega\rightarrow \mathbb{R}^2$ where $\Omega$ is an open subset of $\mathbb{R}^2$ such that $D\subset\Omega$. In the second step, using these decay estimates, an $L^1$-convergence result for bivariate Chebyshev approximant of a function $f$, which satisfies certain regularity conditions on both the variables, is derived.
For this, let us define
\begin{equation}\label{eq:intVitaliTV}
V_{1,1} = \int_{D}\left|f_{xy}(x,y)\right|\omega(x,y)dxdy,
\end{equation}
which is a double Stiltjes integral defined for functions of bounded variation in the sense of Vitali (for details, see \cite{cla-33a}). 
Please observe that the value of $V_{1,1}$ can be infinite depending on the behavior of $x,y=\pm 1$; however, for our results, we are interested only in the case when $V_{k,l}<\infty$ for some integers $k,l\ge 0$.

Let us define
\begin{equation}\label{eq:c_ij^rs}
c_{i,j}^{(r,s)} = \dfrac{4}{\pi^2}\int_{0}^{\pi}\int_{0}^{\pi}f_{x^{r}y^s}(\cos\theta_x,\cos\theta_y)\cos i\theta_x \cos j\theta_y d\theta_x d\theta_y,
\end{equation}
for some integer $r,s\ge0$
and $c_{i,j}^{(0,0)} = c_{i,j}$.

The following theorem is a (modified) generalization of the results obtained by Majidian \cite{has-17a} and Trefethen \cite{tre-13a} for univariate Chebyshev series coefficients. 
\begin{theorem}\label{thm:decaybounds.2dche}
Let $f:D\rightarrow \mathbb{R}$ be a function such that for some integers $k,l\ge 0$, $f_{x^ky^l}$ be of bounded variation in the sense of Vitali on $D$ and $c_{i,j}^{(k+1,l+1)}$ are well defined.
	If 
	$$V_{k,l}[x,y]:= \int_{D}\left|f_{x^{(k+1)}y^{(l+1)}}(x,y)\right|\omega(x,y)dxdy <\infty$$ 
	then for $i\ge k+1$ and $j\ge l+1$, we have
	\begin{equation}\label{eq:bicoeffboundmajidian}
|c_{i,j}| \le \dfrac{4V_{k,l}}{\pi^2}\begin{cases}\Gamma_{0,0}[s](i)\Gamma_{0,0}[r](j), ~~~~~~~\mbox{if }k=2s, l=2r,\\
\Gamma_{0,0}[s](i)\Gamma_{1,-1}[r](j),~~~~~~\mbox{if }k=2s, l=2r+1,\\
\Gamma_{1,-1}[s](i)\Gamma_{0,0}[r](j),~~~~~~\mbox{if }k=2s+1, l=2r,\\
\Gamma_{1,-1}[s](i)\Gamma_{1,-1}[r](j),~~~~~~\mbox{if }k=2s+1, l=2r+1,
\end{cases}
\end{equation}
where $$\Gamma_{\alpha,\beta}[p](\eta)=\dfrac{1}{
\displaystyle \prod_{n=-p}^{p+\alpha} (\eta+2n+\beta) }$$ and 
 $s,r\ge 0$ are integers.
\end{theorem}
\textbf{Proof: } 
Employing integration by parts in \eqref{eq:c_ij^rs} with respect to $\theta_x$ and using $2\sin\theta_x\sin i\theta_x = \cos(i-1)\theta_x-\cos(i+1)\theta_x$ yields
\begin{align}\label{eq:cij^r}
c_{i,j}^{(r,s)}&=\dfrac{1}{2i}\left(c_{i-1,j}^{(r+1,s)}-c_{i+1,j}^{(r+1,s)}\right)
\end{align}
for $r = 0,1,\ldots,k, s = 0,1,\ldots,l+1$ and $i = 1,2,\ldots, j=0,1,2,\ldots$. Similarly, employing integration by parts in \eqref{eq:c_ij^rs} with respect to $\theta_y$  yields
\begin{align}\label{eq:cij^s}
c_{i,j}^{(r,s)}&=\dfrac{1}{2j}\left(c_{i,j-1}^{(r,s+1)}-c_{i,j+1}^{(r,s+1)}\right)
\end{align}
for $r = 0,1,\ldots,k+1, s = 0,1,\ldots,l$ and $i=0,1,2\ldots, j = 1,2,\ldots$. 


To prove the required estimate \eqref{eq:bicoeffboundmajidian}, we first prove the following general inequality 
\begin{equation}\label{eq:bicoeffboundgeneral}
|c_{i,j}^{(k-n,l-m)}| \leq \dfrac{4V_{k,l}}{\pi^2}\begin{cases}\Gamma_{0,0}[s](i)\Gamma_{0,0}[r](j), ~~~~~~~\mbox{if }n=2s, m=2r,\\
\Gamma_{0,0}[s](i)\Gamma_{1,-1}[r](j), ~~~~~~~\mbox{if }n=2s, m=2r+1,\\
\Gamma_{1,-1}[s](i)\Gamma_{0,0}[r](j), ~~~~~~~\mbox{if }n=2s+1, m=2r,\\
\Gamma_{1,-1}[s](i)\Gamma_{1,-1}[r](j), ~~~~~~~\mbox{if }n=2s+1, m=2r+1,
\end{cases}
\end{equation}
for $n = 0,1,\ldots,k, m=0,1,\ldots,l$ and $i\ge n+1, j \geq m+1$. Then $n=k, m=l$ gives the required result. 

Taking $r=k+1$ and $s=l+1$ in \eqref{eq:c_ij^rs}, we get
$$|c_{i,j}^{(k+1,l+1)}| \le \dfrac{4}{\pi^2}\int_{0}^{\pi}\int_{0}^{\pi}|f_{x^{k+1}y^{l+1}}(\cos\theta_x,\cos\theta_y)|d\theta_x d\theta_y = \dfrac{4V_{k,l}}{\pi^2},$$
since $d\theta_x = dx/\sqrt{1-x^2}$ and $d\theta_y = dy/\sqrt{1-y^2}$. Similarly, substituting $r=k,s=l+1$ in \eqref{eq:cij^r} and $r=k+1,s=l$ in \eqref{eq:cij^s}, we get
$$|c_{i,j}^{(k,l+1)}|\le\dfrac{1}{2i}(|c_{i-1,j}^{(k+1,l+1)}|+|c_{i+1,j}^{(k+1,l+1)}|)\le\dfrac{4V_{k,l}}{\pi^2i},~~~~i\ge1,j\ge0$$
and 
$$|c_{i,j}^{(k+1,l)}|\le\dfrac{1}{2j}(|c_{i,j-1}^{(k+1,l+1)}|+|c_{i,j+1}^{(k+1,l+1)}|)\le\dfrac{4V_{k,l}}{\pi^2j},~~~~i\ge0,j\ge1,$$
respectively. We prove \eqref{eq:bicoeffboundgeneral} by double induction on $n$ and $m$. 
\begin{enumerate}
\item For $n=m=0$ case, add \eqref{eq:cij^r} and \eqref{eq:cij^s}, and substitute $r=k, s=l$, we get 
	\begin{align*}
	|c_{i,j}^{(k,l)}|& \le \dfrac{1}{4}\left(\dfrac{1}{i}(|c_{i-1,j}^{(k+1,l)}|+|c_{i+1,j}^{(k+1,l)}|)+\dfrac{1}{j}(|c_{i,j-1}^{(k,l+1)}|+|c_{i,j+1}^{(k,l+1)}|)\right)\\
	& \le\dfrac{1}{4}\left(\dfrac{8V_{k,l}}{\pi^2ij}+\dfrac{8V_{k,l}}{\pi^2ij}\right)=\dfrac{4V_{k,l}}{\pi^2ij},~~~~i,j\ge 1.
	\end{align*}
	\item Let us assume that the inequality \eqref{eq:bicoeffboundgeneral} is true for $n = 2s,m=0$, $i-1\ge 2s+1$ and $j\ge1$. Then for $n=2s+1, s\ge1$ (odd), $m=0$,
	%
	%
	%
	we have
	\begin{align*}
	|c_{i,j}^{(k-2s-1,l)}| & \leq \dfrac{1}{2i}\left(|c_{i-1,j}^{(k-2s,l)}| + |c_{i+1,j}^{(k-2s,l)}|\right)\\ 
&\le\dfrac{1}{2i} \frac{4V_{k,l}}{\pi^2j}\left(\Gamma_{0,-1}[s](i)+\Gamma_{0,1}[s](i)\right)= \dfrac{4V_{k,l}}{\pi^2j}\Gamma_{1,-1}[s](i).
	\end{align*}
	
	\item Assume that the inequality \eqref{eq:bicoeffboundgeneral} holds for $n = 2s+1$ and $m=0$, $i-1\ge 2s+2$ and $j\ge1$. Then for $n=2s+2$ (even) $m=0$, we have (using \eqref{eq:cij^r})
	\begin{align*}
	|c_{i,j}^{(k-2s-2,l)}| & \le \dfrac{1}{2i}\left(|c_{i-1,j}^{(k-2s-1,l)}| + |c_{i+1,j}^{(k-2s-1,l)}|\right)\\ 
	& \leq \dfrac{1}{2i} \frac{4V_{k,l}}{\pi^2j}\left(\Gamma_{1,-2}[s](i)+\Gamma_{1,0}[s](i)\right)= \dfrac{4V_{k,l}}{\pi^2j}\Gamma_{0,-2}[s](i).
	\end{align*}
	\item Similarly, we can prove the inequality \eqref{eq:bicoeffboundgeneral} by fixing $n=0$ and taking $m$ even or odd.
	\item Assume that the inequality \eqref{eq:bicoeffboundgeneral} holds for $n = 2s$ and $m=2r$, $i-1\ge 2s$ and $j-1\ge2r+1$. Then for $n=2s$ and $m=2r+1$ (odd), we have (using \eqref{eq:cij^s})
	\begin{align*}
	|c_{i,j}^{(k-2s,l-2r-1)}| & \leq \dfrac{1}{2j}\left(|c_{i,j-1}^{(k-2s,l-2r)}| + |c_{i,j+1}^{(k-2s,l-2r)}|\right)\\ 
	& \leq \dfrac{1}{2j} \frac{4V_{k,l}}{\pi^2}\Gamma_{0,0}[s](i)\left(\Gamma_{0,-1}[r](j)+\Gamma_{0,1}[r](j)\right)= \dfrac{4V_{k,l}}{\pi^2}\Gamma_{0,0}[s](i)\Gamma_{1,-1}[r](j).
	\end{align*}
	\item Assume that the inequality \eqref{eq:bicoeffboundgeneral} holds for $n = 2s$ and $m=2r+1$, $i-1\ge 2s$ and $j-1\ge2r+2$. Then for $n=2s$ and $m=2r+2, r\ge 1$ (even), we have (using \eqref{eq:cij^s})
	\begin{align*}
|c_{i,j}^{(k-2s,l-2r-2)}| & \leq \dfrac{1}{2j}\left(|c_{i,j-1}^{(k-2s,l-2r-1)}| + |c_{i,j+1}^{(k-2s,l-2r-1)}|\right)\\ 
& \leq \dfrac{1}{2j} \frac{4V_{k,l}}{\pi^2}\Gamma_{0,0}[s](i)\left(\Gamma_{1,-2}[r](j)+\Gamma_{1,0}[r](j)\right)
= \dfrac{4V_{k,l}}{\pi^2}\Gamma_{0,0}[s](i)\Gamma_{0,-2}[r](j).
	\end{align*}
	
\item Similarly by assuming that the inequality \eqref{eq:bicoeffboundgeneral} holds for $n = 2s+1$, we can prove the inequality taking $m$ even or odd. Now the required result is true by induction.
\end{enumerate}
\begin{remark}
In Theorem \ref{thm:decaybounds.2dche}, we obtained element-wise bounds for Chebyshev coefficients $c_{i,j}$. The bounds essentially are of the form $\Gamma[s](i)\Gamma[r](j)$, where $\Gamma[p](\eta)$ converges to zero as $\eta\rightarrow \infty$. We can therefore approximate the coefficient tensor by setting all the entries to zero for which the bounds are lower than a threshold $\epsilon>0$ and subsequently, obtain a low-rank approximation (for detailed survey, see \cite{kis_sch-17a}) for bivariate functions. The extension of the above result in a 3D (three dimensions) setting is also possible, and in that case, the singular values of the matrizations of the coefficient tensors can be used to bound the error for Tucker approximations\footnote{low-rank approximation via Tucker decomposition is used to approximate multivariate functions in the context of Chebfun (for trivariate functions, see \cite{dol-kres-str_21a,has_tre-17a}) and for multivariate functions see, \cite{sai_min_kil-21a}.}.
\end{remark}
\begin{corollary}
Assume the hypothesis of Theorem \ref{thm:decaybounds.2dche} for some positive integers $k$ and $l$. Then for given integers $d_x\ge k$ and $d_y\ge l$, we have
\begin{equation}\label{eq:ApprErrExct}
\|f-{\mathsf{C}}_{d_x,d_y}[f]\|_1 \le \dfrac{4V_{k,l}}{kl\pi^2}\begin{cases}\Pi_{1}[s](d_x)\Pi_{1}[r](d_y),~~~~~~~~~\mbox{if }k=2s, l=2r,\\
\Pi_{1}[s](d_x)\Pi_{0}[r](d_y),
~~~~~~~~~\mbox{if }k=2s, l=2r+1,\\
\Pi_{0}[s](d_x)\Pi_{1}[r](d_y),~~~~~~~~~\mbox{if }k=2s+1, l=2r,\\
\Pi_{0}[s](d_x)\Pi_{0}[r](d_y),~~~~~~~~~\mbox{if }k=2s+1, l=2r+1,
\end{cases}
\end{equation}
for some integer $s,r\ge 0$, where
$$\Pi_{\alpha}[p](n^*)=\dfrac{1}{\displaystyle\prod_{m=-p}^{p-\alpha}(n^*+2m+\alpha)}+\dfrac{1}{\displaystyle\prod_{m=-p}^{p-\alpha}(n^*+2m+\alpha+1)}.$$
\end{corollary}
\textbf{Proof: } For $(x,y)\in D$
\begin{equation}
\|f - {\mathsf{C}}_{d_x,d_y}[f]\|_1 
= \int_D \left|f(x,y) - {\mathsf{C}}_{d_x,d_y}[f](x,y)\right|dxdy \le 4 \sum_{i = d_x+1}^{\infty} \sum_{j = d_y+1}^{\infty} \left|c_{i,j}\right|,\nonumber
\end{equation}
since $\int_{D}dxdy=4$. Using the decay bounds from Theorem \ref{thm:decaybounds.2dche} and then the telescopic property of the resulting series in the above expression, we get the required results.
\begin{theorem}\label{thm:decayboundsforx.2dche}
Let $f:\Omega\rightarrow \mathbb{R}$ be a function such that for some integers $k,l\ge 0$,  $f_{x^k}(x,y_0)$ and $f_{y^l}(x_0,y)$ be of bounded variation on $D$ for $y_0\in [-1,1]$ and $x_0\in[-1,1]$, respectively.
\begin{enumerate}
\item If $V_{k}[x]:=\int_D|f_{x^{(k+1)}}|\omega(x,y)dxdy <\infty,$ 
then for each nonnegative integer $j\in \mathbb{Z}$ and $i\ge k+1$, we have
\begin{equation}
|c_{i,j}| \le \dfrac{4V_k}{\pi^2}\begin{cases}
\Gamma_{0,0}[s](i), ~~~~~~~~~\mbox{if }k=2s,\\
\Gamma_{1,-1}[s](i),~~~~~~~~~\mbox{if }k=2s+1.
\end{cases}
\end{equation}
\item If 
$V_{l}[y]:=\int_D|f_{y^{(l+1)}}|\omega(x,y)dxdy <\infty$ 
then for each nonnegative integer $i\in \mathbb{Z}$ and $j\ge l+1$, we have
\begin{equation}
|c_{i,j}| \le \dfrac{4V_l}{\pi^2}\begin{cases}
\Gamma_{0,0}[r](j), ~~~~~~~~~\mbox{if }l=2r,\\
\Gamma_{1,-1}[r](j),~~~~~~~~~~\mbox{if }l=2r+1.
\end{cases}
\end{equation}
\end{enumerate}	
\end{theorem}
\textbf{Proof: } 
The required result can be obtained by following the same steps as in the univariate case presented in \cite{has-17a} by treating $f_{x^k}$ (for the first case) and $f_{y^l}$(for the second case) as a function of $x$ and $y$, respectively.

\subsubsection*{Relation between $\tilde{c}$ and $c$}
For $(x,y)\in D$, assume the representation
\begin{equation}\label{eq:xseries}
f(x,y) \sim \sum_{j=0}^{\infty}{\vphantom{\sum}}'g_j(x)T_j(y),
\mbox{ where} ~~~
g_j(x) = \sum_{i=0}^{\infty}{\vphantom{\sum}}'c_{i,j}T_i(x).
\end{equation}
the prime ($'$) on the summation denotes that the first term is halved. 
%
%
From \eqref{eq:bivarChebcoeffapprx} we can write
\begin{equation*}
\tilde{c}_{i,j} = \dfrac{2}{n_x}\sum_{l_x=1}^{n_x}\tilde{g}_j(x_{l_x})T_i(x_{l_x}), \mbox{ where }~~\tilde{g}_j(x):= \dfrac{2}{n_y}\sum_{l_y=1}^{n_y}f(x,y_{l_y})T_j(y_{l_y})\approx {g}_j(x).
\end{equation*}
For every $x\in [-1,1]$, we have (for details, see \cite[p. 149]{riv-74a})
\begin{equation}\label{eq:ycoeffrelation}
\tilde{g}_j(x) = g_j(x) + \sum_{k_y=1}^{\infty}(-1)^{k_y}\left(g_{2k_yn_y-j}(x)+g_{2k_yn_y+j}(x)\right).
\end{equation}
On substituting the value of $g_j(x)$ from \eqref{eq:xseries} to \eqref{eq:ycoeffrelation} and resultant value of $\tilde{g}_j(x)$ in $\tilde{c}_{i,j}$, we get
\begin{multline}\label{ExactApproximateRelation2D.eq}
\tilde{c}_{i,j} = c_{i,j}+\sum_{k_x=1}^{\infty}(-1)^{k_x}\left(c_{2k_xn_x-i,j}+c_{2k_xn_x+i,j}\right)+\sum_{k_y=1}^{\infty}(-1)^{k_y}\left(c_{i,2k_yn_y-j}+c_{i,2k_yn_y+j}\right)\\
+\sum_{k_x=1}^{\infty}\sum_{k_y=1}^{\infty}(-1)^{k_x+k_y}\left(c_{2k_xn_x-i,2k_yn_y-j}+c_{2k_xn_x-i,2k_yn_y+j}+c_{2k_xn_x+i,2k_yn_y-j}+c_{2k_xn_x+i,2k_yn_y+j}\right).
\end{multline}
\begin{corollary}\label{thm:errest.2dche}
	Assume the hypothesis of Theorem \ref{thm:decaybounds.2dche} and \ref{thm:decayboundsforx.2dche} for some integers $l,k\ge 1$,
	Then, for given integers $d_x,d_y,n_x$, and $n_y$ such that $n_x-1\geq k,n_y-1\ge l$ and $k\leq d_x <n_x$, $l\leq d_y <n_y$, we have
	\begin{equation}
	\|f-\tilde{\mathsf{C}}_{d_x,d_y}[f]\|_1 \le\dfrac{8V^*}{\pi^2}\left(\dfrac{4}{kl(d_x-k+1)^{k}(d_y-l+1)^{l}}+\dfrac{2(d_y+1)}{k(d_x-k+1)^{k}}+\dfrac{2(d_x+1)}{l(d_y-l+1)^{l}}\right)
	\end{equation}
where $V^* = \max\{V_{k,l},V_k,V_l\}.$
\end{corollary}

\textbf{Proof:} For any $(x,y) \in D$ we have 
\begin{align}
\|f - \tilde{\mathsf{C}}_{d_x,d_y}[f]\|_1 
&\leq \int_D\left |f(x,y) - \mathsf{C}_{d_x,d_y}[f](x,y)\right|dxdy + \int_D \left|{\mathsf{C}}_{d_x,d_y}[f](x,y) - \tilde{\mathsf{C}}_{d_x,d_y}[f](x,y)\right|dxdy,\nonumber\\
&\leq 4 \left(\sum_{i = d_x+1}^{\infty} \sum_{j = d_y+1}^{\infty} \left|c_{i,j}\right|+\sum_{i = 0}^{d_x}{\vphantom{\sum}}'\sum_{j = 0}^{d_y}{\vphantom{\sum}}' \left|c_{i,j}-\tilde{c}_{i,j}\right| \right),\nonumber
\end{align}
since $\int_{D}dxdy=4$.
From the relation \eqref{ExactApproximateRelation2D.eq}, we get
\begin{multline*}
\sum_{i = 0}^{d_x}{\vphantom{\sum}}'\sum_{j = 0}^{d_y}{\vphantom{\sum}}'\left|c_{i,j}-\tilde{c}_{i,j}\right|\le
\sum_{j = 0}^{d_y}{\vphantom{\sum}}'\sum_{k_x=1}^{\infty}\sum_{i=2k_xn_x-d_x}^{2k_xn_x+d_x}|c_{i,j}|+\sum_{i = 0}^{d_x}{\vphantom{\sum}}'\sum_{k_y=1}^{\infty}\sum_{j=2k_yn_y-d_y}^{2k_yn_y+d_y}|c_{i,j}|
+\sum_{k_x=1}^{\infty}\sum_{k_y=1}^{\infty}\sum_{i=2k_xn_x-d_x}^{2k_xn_x+d_x}\sum_{j=2k_yn_y-d_y}^{2k_yn_y+d_y}|c_{i,j}|.
\end{multline*}
For $d_x=n_x-l_x,l_x=1,2,\ldots,n_x-k$ and $d_y=n_y-l_y,l_y=1,2,\ldots,n_y-l$, we can write
\begin{align}\label{eq:ConvIntStep}
\|f - \tilde{\mathsf{C}}_{d_x,d_y}[f]\|_1 \le & 4\left(\sum_{i = n_x-l_x+1}^{\infty} \sum_{j = n_y-l_y+1}^{\infty} \left|c_{i,j}\right|+\sum_{k_x=1}^{\infty}\sum_{k_y=1}^{\infty}\sum_{i=(2k_x-1)n_x+l_x}^{(2k_x+1)n_x-l_x}\sum_{j=(2k_y-1)n_y+l_y}^{(2k_y+1)n_y-l_y}|c_{i,j}|\right.\nonumber\\
&\left.+\sum_{j = 0}^{d_y}{\vphantom{\sum}}'\sum_{k_x=1}^{\infty}\sum_{i=(2k_x-1)n_x+l_x}^{(2k_x+1)n_x-l_x}|c_{i,j}|+\sum_{i = 0}^{d_x}{\vphantom{\sum}}'\sum_{k_y=1}^{\infty}\sum_{j=(2k_y-1)n_y+l_y}^{(2k_y+1)n_y-l_y}|c_{i,j}|\right)\nonumber\\
\le & 4\left(2\sum_{i = n_x-l_x+1}^{\infty} \sum_{j = n_y-l_y+1}^{\infty} \left|c_{i,j}\right|+\sum_{j = 0}^{n_y-l_y}{\vphantom{\sum}}'\sum_{i = n_x-l_x+1}^{\infty} \left|c_{i,j}\right|+\sum_{i = 0}^{n_x-l_x}{\vphantom{\sum}}' \sum_{j = n_y-l_y+1}^{\infty} \left|c_{i,j}\right|\right).
\end{align}
From the hypothesis, we have $d_x(=n_x-l_x)\ge k$ and $d_y(=n_y-l_y)\ge l$ hence, we can use the decay bounds from Theorem \ref{thm:decaybounds.2dche} and \ref{thm:decayboundsforx.2dche} in the above expression. Using the decay bounds and the telescopic property of the resulting series, we get 
\begin{equation}\label{eq:turnerrbichemixed.2dche}
\sum_{i = n_x-l_x+1}^{\infty} \sum_{j = n_y-l_y+1}^{\infty} \left|c_{i,j}\right| \le \dfrac{4V_{k,l}}{4kl\pi^2}\begin{cases}\Pi_{1}[s](n_x-l_x)\Pi_{1}[r](n_y-l_y),\\~~~~~~~~~\mbox{if }k=2s, l=2r,\\
\Pi_{1}[s](n_x-l_x)\Pi_{0}[r](n_y-l_y),\\
~~~~~~~~~\mbox{if }k=2s, l=2r+1,\\
\Pi_{0}[s](n_x-l_x)\Pi_{1}[r](n_y-l_y),\\~~~~~~~~~\mbox{if }k=2s+1, l=2r,\\
\Pi_{0}[s](n_x-l_x)\Pi_{0}[r](n_y-l_y),\\ ~~~~~~~~~\mbox{if }k=2s+1, l=2r+1,
\end{cases}
\end{equation}
\begin{equation}\label{eq:turnerrbichex.2dche}
\sum_{j = 0}^{n_y-l_y}{\vphantom{\sum}}'\sum_{i = n_x-l_x+1}^{\infty} \left|c_{i,j}\right|\le (n_y-l_y+1)\dfrac{4V_k}{2k\pi^2}\begin{cases}\Pi_{1}[s](n_x-l_x),~~~~~~~~~\mbox{if }k=2s,\\
\Pi_{0}[s](n_x-l_x),~~~~~~~~~\mbox{if }k=2s+1,
\end{cases}
\end{equation}
and
\begin{equation}\label{eq:turnerrbichey.2dche}
\sum_{i = 0}^{n_x-l_x}{\vphantom{\sum}}'\sum_{j = n_y-l_y+1}^{\infty} \left|c_{i,j}\right|\le (n_x-l_x+1)\dfrac{4V_l}{2l\pi^2}\begin{cases}\Pi_{1}[r](n_y-l_y),~~~~~~~~~\mbox{if }l=2r,\\
\Pi_{0}[r](n_y-l_y),~~~~~~~~~\mbox{if }l=2r+1,
\end{cases}
\end{equation}
where
$$\Pi_{\alpha}[p](n^*)=\dfrac{1}{\displaystyle\prod_{m=-p}^{p-\alpha}(n^*+2m+\alpha)}+\dfrac{1}{\displaystyle\prod_{m=-p}^{p-\alpha}(n^*+2m+\alpha+1)}\le \dfrac{2}{(n^*-2p+\alpha)^{2p-\alpha+1}}.$$
Substituting the estimates from \eqref{eq:turnerrbichemixed.2dche}, \eqref{eq:turnerrbichex.2dche} and \eqref{eq:turnerrbichey.2dche} into \eqref{eq:ConvIntStep} and using the above estimate for $\Pi_{\alpha}[p](n^*)$ with appropreate values of $p,n^*$ and $\alpha$, we get the required result.
\section*{Conclusion} We obtained decay estimates and $L^1$-convergence results for bivariate Chebyshev approximations for functions of limited regularity. In general, low-rank approximations of the coefficient tensor are used to approximate 2D functions (see \cite{tow_tre-13a}). However, we are inclined towards the nonlinear approximation methods to approximate functions with discontinuities and using these decay bounds, one can obtain convergence results for Chebyshev-based nonlinear approximation methods, for instance, the bivariate Pad\'e-Chebyshev method. 
%
%
\bibliography{bibfile_thesis}

\begin{thebibliography}{10}

\bibitem{ada_cla-34a}
C~R Adams and J~A Clarkson.
\newblock Properties of functions f (x, y) of bounded variation.
\newblock {\em Transactions of the American Mathematical Society},
  36(4):711--730, 1934.

\bibitem{cla-33a}
J.~A. Clarkson.
\newblock On double riemann-stieltjes integrals.
\newblock {\em Bulletin of the American Mathematical Society}, 39(12):929--936,
  1933.

\bibitem{cla_ada-33a}
J.~A. Clarkson and C.~R. Adams.
\newblock On definitions of bounded variation for functions of two variables.
\newblock {\em Transactions of the American Mathematical Society},
  35(4):824--854, 1933.

\bibitem{dol-kres-str_21a}
S.~Dolgov, D.~Kressner, and C.~Strössner.
\newblock Functional tucker approximation using chebyshev interpolation.
\newblock {\em SIAM Journal on Scientific Computing}, 43(3):A2190--A2210, 2021.

\bibitem{dri_hal_tre-14a}
T.~A. Driscoll, N.~Hale, and L.~N. Trefethen.
\newblock Chebfun guide.
\newblock {\em Pafnuty Publications, Oxford, UK,}, 2014.

\bibitem{gol_kal_lut-13a}
R.~V. Golovanov, N.~N. Kalitkin, and K.~I. Lutskiy.
\newblock Odd extension for fourier approximation of nonperiodic functions.
\newblock {\em Mathematical Models and Computer Simulations}, 5(6):595--606,
  2013.

\bibitem{has_tre-17a}
B.~Hashemi and L.~N. Trefethen.
\newblock Chebfun in three dimensions.
\newblock {\em SIAM Journal on Scientific Computing}, 39(5):C341--C363, 2017.

\bibitem{kis_sch-17a}
N~Kishore~K. and J.~Schneider.
\newblock Literature survey on low rank approximation of matrices.
\newblock {\em Linear and Multilinear Algebra}, 65(11):2212--2244, 2017.

\bibitem{has-17a}
Hassan M.
\newblock On the decay rate of chebyshev coefficients.
\newblock {\em Applied Numerical Mathematics}, 113:44 -- 53, 2017.

\bibitem{mas-han-03a}
J.~C. Mason and D.~C. Handscomb.
\newblock {\em Chebyshev polynomials}.
\newblock Chapman \& Hall/CRC, Boca Raton, FL, 2003.

\bibitem{pac_pla_tre-10a}
R.~Pach\'{o}n, R.~B. Platte, and L.~N. Trefethen.
\newblock Piecewise-smooth chebfuns.
\newblock {\em IMA J. Numer. Anal.}, 30(4):898--916, 2010.

\bibitem{riv-74a}
T.~J. Rivlin.
\newblock {\em The {C}hebyshev polynomials}.
\newblock Wiley-Interscience [John Wiley \& Sons], New York-London-Sydney,
  1974.
\newblock Pure and Applied Mathematics.

\bibitem{sai_min_kil-21a}
A.~K Saibaba, R.~Minster, and M.~E Kilmer.
\newblock Efficient randomized tensor-based algorithms for function
  approximation and low-rank kernel interactions.
\newblock {\em arXiv preprint arXiv:2107.13107}, 2021.

\bibitem{tow_tre-13a}
A.~Townsend and L.~N. Trefethen.
\newblock An extension of chebfun to two dimensions.
\newblock {\em SIAM Journal on Scientific Computing}, 35(6):C495--C518, 2013.

\bibitem{tre-08a}
L.~N. Trefethen.
\newblock Is {G}auss quadrature better than {C}lenshaw-{C}urtis?
\newblock {\em SIAM Rev.}, 50(1):67--87, 2008.

\bibitem{tre-13a}
L.~N. Trefethen.
\newblock {\em Approximation theory and approximation practice}.
\newblock Society for Industrial and Applied Mathematics (SIAM), Philadelphia,
  PA, 2013.

\bibitem{xia-18a}
S.~Xiang.
\newblock On the optimal convergence rates of {C}hebyshev interpolations for
  functions of limited regularity.
\newblock {\em Appl. Math. Lett.}, 84:1--7, 2018.

\bibitem{xia_che_wan-10a}
S.~Xiang, X.~Chen, and H.~Wang.
\newblock Error bounds for approximation in chebyshev points.
\newblock {\em Numerische Mathematik}, 116(3):463--491, Sep 2010.

\end{thebibliography}
\end{document}